\documentclass[11pt]{article}

\usepackage{amssymb,amsmath,amsfonts,amsthm}
\usepackage{graphicx,psfrag}

\newcommand{\N}{{\mathbb N}}
\newcommand{\R}{{\mathbb R}}
\newcommand{\SP}{{\mathbb S}}
\newcommand{\W}{{\mathbb W}}

\newtheorem{thh}{Theorem}
\newtheorem{lemm}{Lemma}
\newtheorem{prop}{Proposition}

\newtheorem{remarque}{Remark}

\begin{document}

\title{Shock Profiles for Non Equilibrium Radiating Gases}

\author{Chunjin {\sc Lin}$^\dag$, Jean-Fran\c{c}ois {\sc Coulombel}$^{\dag \ddag}$, 
Thierry {\sc Goudon}$^{\dag \ddag}$\\
\\
{\small $\dag$ Team SIMPAF--INRIA Futurs \& Universit\'e Lille 1, 
Laboratoire Paul Painlev\'e, UMR CNRS 8524}\\
{\small Cit\'e scientifique, 59655 VILLENEUVE D'ASCQ Cedex, France}\\
{\small $\ddag$ CNRS}\\
{\small E-mails: {\tt chunjin.lin@math.univ-lille1.fr}, 
{\tt jfcoulom@math.univ-lille1.fr}},\\
{\small \tt thierry.goudon@math.univ-lille1.fr}}
\date{\today}
\maketitle

\begin{abstract}
We study a model of radiating gases that describes the interaction of an 
inviscid gas with photons. We show the existence of smooth traveling waves 
called 'shock profiles', when the strength of the shock is small. Moreover, 
we prove that the regularity of the traveling wave increases when the 
strength of the shock tends to zero.
\end{abstract}

\section{Introduction and main results}

We are interested in a system of PDEs describing astrophysical flows, where 
a gas interacts with radiation through energy exchanges. Similar questions 
arise in the modeling of reentry problems, or high temperature combustion 
phenomena. The gas is described by its density $\rho >0$, its bulk velocity 
$u \in \R$, and its specific total energy $E=e+u^2/2$, where $e$ stands for 
the specific internal energy. (Our analysis is restricted to a one-dimensional 
framework, but this is not a loss of generality, as shown below.) We consider 
a situation where the gas is not in thermodynamical equilibrium with the 
radiations, which are thus described by their own energy $n$. The evolution 
of the gas flow is governed by the system:
\begin{equation}
\label{eulerevol}
\begin{cases}
\partial_t \rho +\partial_x (\rho \, u)=0 \, ,& \\
\partial_t (\rho \, u) +\partial_x (\rho \, u^2+P)=0\, ,& \\
\partial_t (\rho \, E) +\partial_x (\rho \, E \, u+P \, u)=n-\theta^4 \, ,& 
\end{cases}
\end{equation}
where the right-hand side in the last equation accounts for energy exchanges 
with the radiations, $P$ being the pressure of the gas, and $\theta$ its 
temperature. Throughout the paper, we always assume that the gas obeys the 
perfect gas pressure law:
\begin{equation}
\label{pressure}
P=R \, \rho \, \theta =(\gamma-1) \, \rho \, e \, ,
\end{equation}
where $R$ is the perfect gas constant, and $\gamma>1$ is the ratio of the 
specific heats at constant pressure, and volume. This assumption yields many 
algebraic simplifications, but we believe that our results still hold for a 
general pressure law satisfying the usual requirements of thermodynamics. 
System \eqref{eulerevol} is completed by considering that radiations are 
described by a stationary diffusion regime that reads:
\begin{equation}
\label{diffn}
-\partial_{xx} n =\theta^4 -n \, .
\end{equation}
We detail in Appendix \ref{model} how the system \eqref{eulerevol}, \eqref{diffn} 
can be formally derived by asymptotics arguments, starting from a more complete 
system involving a kinetic equation for the specific intensity of radiation.

As a matter of fact, the operator $(1-\partial_{xx})$ can be explicitly 
inverted, and \eqref{diffn} can be recast as a convolution:
\begin{equation}
\label{diffn2}
n(t,x) =\dfrac{1}{2} \, \displaystyle \int_{\R} 
{\rm e}^{-|x-y|} \, \theta (t,y)^4 \, dy \, .
\end{equation}
Let us introduce the quantity:
\begin{equation}
\label{defq}
q(t,x) :=-\partial_x n \, (t,x)=\dfrac{1}{2} \displaystyle \int_{\R} 
{\rm e}^{-|x-y|} \, {\rm sgn}(x-y) \, \theta (t,y)^4 \, dy \, ,
\end{equation}
where sgn is the sign function:
\begin{equation*}
{\rm sgn} (x)=\begin{cases}
1 &\text{if $x>0$,} \\
0 &\text{if $x=0$,} \\
-1 &\text{if $x<0$.}
\end{cases}
\end{equation*}
The quantity $q$ can be interpreted as the radiative heat flux. Then, we 
can rewrite \eqref{eulerevol}, \eqref{diffn} as follows:
\begin{equation}
\label{Euler2}
\begin{cases}
\partial_t \rho +\partial_x (\rho \, u) =0\, ,& \\
\partial_t (\rho \, u) +\partial_x (\rho \, u^2+P)=0 \, ,& \\
\partial_t (\rho \, E) +\partial_x (\rho \, E \, u +P \, u +q) =0 \, ,& 
\end{cases}
\end{equation}
with $q$ given by \eqref{defq}. Recall that $E=e+u^2/2$, and $P$ is given 
by \eqref{pressure}.

In this paper, we address the question of the influence of the energy exchanges 
on the structure of shock waves. More precisely, let us consider given states at 
infinity $(\rho_\pm, u_\pm,e_\pm)$, and let us asume that:
\begin{equation}
\label{shock}
(\rho,u,e)(t,x) =\begin{cases}
(\rho_-,u_-,e_-) &\text{if $x<\sigma \, t$,} \\
(\rho_+,u_+,e_+) &\text{if $x>\sigma \, t$,}
\end{cases}
\end{equation}
is a shock wave, with speed $\sigma$, solution to the standard Euler 
equations (that is, system \eqref{Euler2} with $q \equiv 0$). We refer to 
\cite{lax,serrelivre,smoller} for a detailed study of shock waves for the 
Euler equations. The question we ask is the following: does there exist a 
traveling wave $(\rho,u,e)(x-\sigma t)$ solution to \eqref{Euler2}, with $q$ 
given by \eqref{defq}, that satisfies the asymptotic conditions:
\begin{equation}
\label{asym-cond}
\lim_{\xi \rightarrow \pm\infty} (\rho,u,e)(\xi) =(\rho_\pm,u_\pm,e_\pm) \, .
\end{equation}
In other words, we are concerned with the existence of a shock profile, and a 
natural expectation (at least for shocks of small amplitude) is that the step 
shock \eqref{shock} is smoothed into a continuous profile, due to the dissipation 
introduced by \eqref{diffn}. The analogous problem for the compressible Navier-Stokes 
system has been treated a long time ago, see \cite{gilbarg}, without any smallness 
assumption on the shock wave. Concerning radiative transfer, a formal analysis 
of shock profiles has been performed in \cite{HB}, together with rough numerical 
simulations. (We refer also to \cite{ZR,MM} for the physical background.) The 
main purpose of this work is to make the analysis of \cite{HB} rigorous. Since 
we are only concerned in this paper with the existence of shock profiles, and not 
with their stability, the problem is purely one-dimensional (due to the Galilean 
invariance of the Euler equations). This is why we have directly restricted to 
the one-dimensional case. However, the formal derivation of Appendix \ref{model} 
is made in several space dimensions.

Before stating our main results, let us mention that a simplified version 
of \eqref{Euler2}, \eqref{defq} has been introduced, and studied in \cite{SchoTad} and later in
\cite{kawanishi,kawanishi2}. This 'baby-model' consists in a Burgers 
type equation:
\begin{equation*}
\partial_t u + \partial_x  \Big( \dfrac{u^2}{2} \Big)=-\partial_x q \, ,
\end{equation*}
coupled to the diffusion equation:
\begin{equation*}
-\partial_{xx} q+q=-\partial_x u \, .
\end{equation*}
These two equations can be seen as a scalar version of \eqref{Euler2}, 
\eqref{defq} since they can be recast as:
\begin{equation}
\label{Kburgers}
\partial_t u +\partial_x \Big( \dfrac{u^2}{2} \Big) =Ku-u \, ,
\end{equation} 
where $K$ is the integral operator already arising in \eqref{defq}:
\[
Ku (t,x)=\dfrac{1}{2} \displaystyle \int_{\R} {\rm e}^{-|x-y|} \, u(t,y) \, dy \, .
\]
The thorough study of \eqref{Kburgers} has motivated a lot of works; we mention 
in particular \cite{nishi,marcati,tadmor,serre}. Clearly \eqref{Kburgers} can be 
seen as a prototype for discussing \eqref{Euler2}, \eqref{defq}; nevertheless, 
replacing \eqref{Euler2}, and \eqref{defq} by \eqref{Kburgers} has two important 
consequences: the equation becomes scalar, and the 'diffusion' $K-1$ applies to 
the unique unknown (while in \eqref{Euler2}, the 'diffusion' appears through the 
radiative heat flux $q$ only in the third equation). Our work is a first attempt 
to extend the known results for \eqref{Kburgers} to the more physical model 
\eqref{Euler2}, \eqref{defq}.

Let us now state our main results. The first result deals with the existence 
of smooth shock profiles when the strength of the shock is small:

\begin{thh}
\label{exis}
Let $\gamma$ satisfy
\begin{equation*}
1< \gamma <\dfrac{\sqrt{7}+1}{\sqrt{7}-1} \simeq 2.215 \, ,
\end{equation*}
and let $(\rho_-,u_-,e_-)$ be fixed. Then there exists a positive constant 
$\delta$ (that depends on $(\rho_-,u_-,e_-)$, and $\gamma$) such that, for 
all state $(\rho_+,u_+,e_+)$ verifying:
\begin{itemize}
\item $\|(\rho_+,u_+,e_+)-(\rho_-,u_-,e_-)\| \le \delta$,

\item the function \eqref{shock} is a shock wave, with speed $\sigma$, for 
the (standard) Euler equations,
\end{itemize}
then there exists a $C^2$ traveling wave $(\rho,u,e)(x-\sigma t)$ solution 
to \eqref{Euler2}, \eqref{defq}, \eqref{asym-cond}.
\end{thh}

As in the study of the 'baby-model' \eqref{Kburgers}, the existence of 
a smooth shock profile is linked to a smallness assumption on the shock 
strength, see \cite{kawanishi}. Here the smallness parameter $\delta$ 
may depend on the state $(\rho_-,u_-,e_-)$, while for \eqref{Kburgers}, 
the smallness parameter is uniform (and even explicit!).

The restriction on the adiabatic constant $\gamma$ might be unnecessary, 
but it simplifies the proof, and it covers the main physical cases 
$1<\gamma \le 2$.

Our second result is also in the spirit of \cite{kawanishi}, and deals with 
the smoothness of the shock profile constructed in the previous Theorem:

\begin{thh}
\label{smooth}
Let $\gamma$ satisfy
\begin{equation*}
1< \gamma < \dfrac{\sqrt{7}+1}{\sqrt{7}-1} \simeq 2.215 \, ,
\end{equation*}
and let $(\rho_-,u_-,e_-)$ be fixed. Then there exists a decreasing 
sequence of positive numbers $(\delta_n)_{n \in \N}$ (the sequence depends 
on $(\rho_-,u_-,e_-)$, and $\gamma$) such that, for all $n\in \N$, and for 
all state $(\rho_+,u_+,e_+)$ verifying:
\begin{itemize}
\item $\|(\rho_+,u_+,e_+)-(\rho_-,u_-,e_-)\| \le \delta_n$,

\item the function \eqref{shock} is a shock wave, with speed $\sigma$, for 
the (standard) Euler equations,
\end{itemize}
then there exists a $C^{n+2}$ traveling wave $(\rho,u,e)(x-\sigma t)$ solution 
to \eqref{Euler2}, \eqref{defq}, \eqref{asym-cond}.
\end{thh}

To a large extent, our analysis follows the arguments of \cite{HB}, \cite{SchoTad} and 
\cite{kawanishi}. The proof of Theorem \ref{exis} is presented in Section 
\ref{Proof}, while Section \ref{moresmooth} is devoted to the proof of 
Theorem \ref{smooth}. The investigation of strong shocks, as well as 
stability issues will be addressed in a forthcoming work.

\section{Existence of smooth shock profiles}
\label{Proof}

In this section, we prove Theorem \ref{exis}. We first recall some basic facts 
on shock waves for the Euler equations. Then, we make some transformations on 
the traveling wave equation. Eventually, we prove Theorem \ref{exis} by using 
an auxiliary system of Ordinary Differential Equations, that is introduced and 
studied in the last paragraph of this section.

\subsection{Shock wave solutions to the Euler equations}

In this paragraph, we recall some basic facts about the (entropic) shock 
wave solutions to the Euler equations:
\begin{equation*}
\begin{cases}
\partial_t \rho +\partial_x (\rho \, u) =0 \, ,\\
\partial_t (\rho \, u) +\partial_x (\rho \, u^2 +P) =0 \, ,\\
\partial_t (\rho \, E) +\partial_x (\rho \, E \, u+ P \, u) = 0 \, ,
\end{cases}
\end{equation*}
where $P$, and $E$ are given as in the introduction. We refer to 
\cite{lax,serrelivre,smoller} for all the details, and omit the calculations. 
In all what follows, we only consider shock waves that satisfy Lax shock 
inequalities. We shall thus speak of $1$-shock waves, or $3$-shock waves.

We consider a fixed 'left' state $(\rho_-,u_-,e_-)$. Then the 'right' states 
$(\rho_+,u_+,e_+)$ such that $(\rho_\pm,u_\pm,e_\pm)$ define a $1$-shock 
wave, with speed $\sigma$, is a half-curve initiating at $(\rho_-,u_-,e_-)$. 
Introducing the notation $v_\pm=u_\pm -\sigma$, the Rankine-Hugoniot jump 
conditions can be rewritten as:
\begin{equation}
\label{defjC1C2}
\begin{cases}
\rho_+ \, v_+ =\rho_- \, v_- =:j \, ,\\
\rho_+ \, v_+^2 +P_+ =\rho_- \, v_-^2 +P_- =:j \, C_1 \, ,\\
\rho_+ \, v_+ \big( e_+ +\dfrac{v_+^2}{2} \big) +P_+ \, v_+ 
=\rho_- \, v_- \big( e_- +\dfrac{v_-^2}{2} \big) +P_- \, v_- =:j \, C_2 \, .
\end{cases}
\end{equation}
Observe that $v_-$ does not only depend on the 'left' state $(\rho_-,u_-,e_-)$, 
but also on $(\rho_+,u_+,e_+)$, because $v_-$ is defined with the help of the 
shock speed $\sigma$. Consequently, the constants $j$, $C_1$, and $C_2$ depend 
on both $(\rho_-,u_-,e_-)$, and $(\rho_+,u_+,e_+)$.

For $1$-shocks, that is when the inequalities:
\begin{equation*}
u_+ -c_+ < \sigma <u_+ \, ,\quad \sigma < u_- -c_- \, ,
\end{equation*}
are satisfied ($c$ denotes the sound speed), all quantities $j$, $C_1$, and 
$C_2$ are positive. Moreover, when the strength of the shock tends to zero, 
that is when $(\rho_+,u_+,e_+)$ tends to $(\rho_-,u_-,e_-)$, one has
\begin{equation}
\label{asymp}
\begin{pmatrix}
\sigma \\
j \\
C_1 \\
C_2 \end{pmatrix} \longrightarrow \begin{pmatrix}
u_- - c_- \\
\rho_- \, c_- \\
c_- +(\gamma-1) e_-/c_- \\
\gamma \, e_- +c_-^2/2 \end{pmatrix} \, .
\end{equation}
Consequently, when the strength of the shock is small, all quantities 
$j$, $C_1$, $C_2$ are bounded away from zero. Recall also that $1$-shocks 
are compressive, in the sense that $\rho_+ >\rho_-$. This inequality 
immediately implies that $0<v_+ <v_-$. Eventually, the strength of the 
shock tends to zero if, and only if $u_+$ tends to $u_-$ (in that case, 
we also have $\rho_+ \rightarrow \rho_-$, and $e_+ \rightarrow e_-$ 
because of the Rankine-Hugoniot jump conditions).

In all what follows, we limit our discussion to the case of $1$-shocks 
for simplicity, but the extension to $3$-shocks is immediate.

\subsection{Reduction of the traveling wave equation}

In this paragraph, we derive, and transform the equation satisfied by 
traveling wave solutions to \eqref{Euler2}, \eqref{defq}. A traveling 
wave solution to \eqref{Euler2}, \eqref{defq} with speed $\sigma$ is a 
solution $(\rho,u,e)(x-\sigma t)$. For such solutions, the radiative 
heat flux $q$ also depends on the sole variable $x-\sigma t$:
\begin{equation*}
q(x-\sigma t)=\dfrac{1}{2} \displaystyle \int_{\R} 
{\rm e}^{-|x-\sigma t-y|} \, {\rm sgn}(x-\sigma t-y) \, \theta (y)^4 \, dy \, ,
\end{equation*}
and \eqref{Euler2} reads:
\begin{equation*}
\begin{cases}
(\rho\, (u-\sigma))' =0 \, ,\\
(\rho \, u \, (u-\sigma)+(\gamma-1) \, \rho \, e)' =0 \, ,\\
(\rho \, (e+\frac{u^2}{2}) \, (u-\sigma)+(\gamma-1) \, \rho \, e \, u+q)' =0 \, ,
\end{cases}
\end{equation*}
where $'$ denotes differentiation with respect to the variable 
$\xi=x-\sigma t$. Introducing the new unknown $v=u-\sigma$, the 
above system is easily seen to be equivalent to:
\begin{equation}
\label{Euler-xi}
\begin{cases}
(\rho \, v)'=0 \, ,\\
(\rho \, v^2+(\gamma-1)\, \rho \, e)' =0 \, ,\\
(\rho \, v \, (e+\frac{v^2}{2})+(\gamma-1) \, \rho \, v \, e+q)' =0 \, .
\end{cases}
\end{equation}
Since we are looking for a shock profile, the traveling wave solution 
should also satisfy:
\begin{equation}
\label{asym-cond2}
\lim_{\xi\to \pm\infty} (\rho,v,e) =(\rho_{\pm},v_{\pm},e_{\pm}) \, ,
\end{equation}
where $v_{\pm}=u_{\pm}-\sigma$, and $(\rho_\pm,u_\pm,e_\pm)$ defines a 
$1$-shock wave with speed $\sigma$ for the Euler equations. Notice that 
the quantity $a=|u_+-u_-|/2$, that measures the strength of the shock, 
is invariant with respect to our change of velocity, that is 
$a=|u_+-u_-|/2=|v_+-v_-|/2$. Recall also that for $1$-shocks, 
there holds $v_->v_+>0$.

Observing that we have:
\begin{equation*}
q(\xi)=\dfrac{1}{2} \displaystyle \int_0^{+\infty} {\rm e}^{-y} \, 
\big( \theta (\xi-y)^4 -\theta (\xi+y)^4 \big) \, dy \, ,
\end{equation*}
we conclude that $q$ tends to zero at $\pm \infty$ by Lebesgue's Theorem 
(because $\theta$ is necessarily bounded since it has finite limits at 
$\pm \infty$). We can thus integrate the system 
\eqref{Euler-xi}-\eqref{asym-cond2} once, and \eqref{Euler-xi} reads 
equivalently:
\begin{equation}
\label{system}
\begin{cases}
(\rho \, v)(\xi)=j \, ,\\
(\rho \, v^2+(\gamma-1) \, \rho \, e)(\xi)=j \, C_1 \, ,\\
(\rho \, v \, (e+\frac{v^2}{2})+(\gamma-1) \, \rho \, v \, e+q)(\xi)=j \, C_2 \, ,
\end{cases}
\end{equation}
where the constants $j,$ $C_1$ and $C_2$ are given by the Rankine-Hugoniot 
conditions \eqref{defjC1C2}. For small shocks, the positive constants $j$, 
$C_1$, $C_2$ have the asymptotic behavior \eqref{asymp}.

From the two first equations of \eqref{system}, we derive the relations:
\begin{equation*}
\rho (\xi)=\dfrac{j}{v(\xi)} \, ,\quad 
e(\xi)=\dfrac{(C_1-v(\xi)) \, v(\xi)}{\gamma-1} \, .
\end{equation*}
The third equation of \eqref{system} thus reduces to:
\begin{equation}
\label{v-q}
v(\xi)^2 -\dfrac{2 \, \gamma \, C_1}{\gamma+1} \, v(\xi) 
+\dfrac{2 \, (\gamma-1) \, C_2}{\gamma+1} 
=\dfrac{2 \, (\gamma-1)}{j \, (\gamma+1)} \, q(\xi) \, .
\end{equation}
Using the equation of state \eqref{pressure}, as well as the second equation 
of \eqref{system}, we get:
\begin{equation*}
\theta (\xi)=\dfrac{(\gamma-1) \, e(\xi)}{R} =\dfrac{(C_1-v(\xi)) \, v(\xi)}{R} \, .
\end{equation*}
Consequently, \eqref{v-q} can be recast as an integral equation with a single 
unknown function $v$:
\begin{equation}
\label{eqv}
v(\xi)^2 -\dfrac{2 \, \gamma \, C_1}{\gamma+1} \, v(\xi) 
+\dfrac{2 \, (\gamma-1) \, C_2}{\gamma+1} 
=\dfrac{(\gamma-1)}{j \, (\gamma+1) \, R^4} \int_\R {\rm e}^{-|\xi-y|} \, 
{\rm sgn} (\xi-y) \, v(y)^4 (C_1-v(y))^4 \, dy \, .
\end{equation}
We are searching for a solution $v$ to \eqref{eqv}, that satisfies the asymptotic 
conditions $v(\xi)\rightarrow v_\pm$, as $\xi \rightarrow \pm \infty$.

\begin{remarque}
If we find a $C^2$ solution $v$ to \eqref{eqv} that does not vanish, and that 
satisfies $v(\xi) \to v_{\pm}$ as $\xi \to \pm \infty$, then we obtain a $C^2$ 
shock profile $(\rho,u,e)$ by simply setting:
\begin{equation*}
\rho(\xi)=\dfrac{j}{v(\xi)} \, ,\quad u(\xi)=v(\xi)+\sigma \, ,\quad 
e(\xi)=\dfrac{(C_1-v(\xi))\, v(\xi)}{\gamma-1} \, .
\end{equation*}
In particular, if $v(\xi) \in [v_+,v_-]$ for all $\xi$, then $v$ does not vanish.
\end{remarque}

\begin{remarque}
Since the heat flux $q$ vanishes at $\pm \infty$, \eqref{eqv} can be also 
rewritten as:
\begin{equation*}
(v(\xi)-v_-)(v(\xi)-v_+) =
\dfrac{(\gamma-1)}{j \, (\gamma+1) \, R^4} \int_\R {\rm e}^{-|\xi-y|} \, 
{\rm sgn} (\xi-y) \, v(y)^4 \, (C_1-v(y))^4 \, dy \, .
\end{equation*}
\end{remarque}

We are going to rewrite \eqref{eqv} as a second order differential equation, 
that will be easier to study than the integral equation \eqref{eqv}. Indeed, 
assuming that $v$ is a $C^2$ function of $\xi$, and differentiating twice 
\eqref{eqv} with respect to $\xi$, we get (see \cite{HB} for the details of 
the computations):
\begin{equation}
\label{eqv1}
(v-\dfrac{\gamma \, C_1}{\gamma+1}) \, v''+(v')^2 
-\dfrac{4 \, (\gamma-1)}{j \, (\gamma+1) \, R^4} \, (C_1-v)^3 \, v^3 
\, (C_1-2v) \, v' -\dfrac{1}{2} \, (v-v_-) \, (v-v_+) =0 \, .
\end{equation}
Conversely, if $v$ is a $C^2$ solution to \eqref{eqv1} that satisfies 
$v(\xi) \to v_{\pm}$ as $\xi \to \pm \infty$, then $v$ is also a solution 
to \eqref{eqv}. If in addition $v$ takes its values in the interval 
$[v_+,v_-]$, then we can construct a $C^2$ shock profile, and thus prove 
Theorem \ref{exis}.

The differential equation \eqref{eqv1} can be simplified by introducing the 
new unknown function $\hat{v}=v-(v_- +v_+)/2$, and by rewriting the second 
order differential equation as a first order system:
\begin{equation}
\label{ode}
\left\{
\begin{array}{lll}
\hat{v}' &=& w \, ,\\
\hat{v} \, w' &=& -w^2-f(\hat{v}) \, w+\dfrac{\hat{v}^2-a^2}{2} \, ,
\end{array}
\right.
\end{equation}
where $f$ is the following polynomial function:
\begin{equation}
\label{fun f}
f(\hat{v})=\dfrac{4 \, (\gamma-1)}{j \, R^4 \, (\gamma+1)} \, 
\left( \dfrac{C_1}{\gamma+1}-\hat{v} \right)^3 \, 
\left( \hat{v}+\dfrac{\gamma\, C_1}{\gamma+1} \right)^3 \, 
\left( 2\hat{v}+\dfrac{(\gamma-1) \, C_1}{\gamma+1} \right) \, .
\end{equation}
We recall that $a=|v_- -v_+|/2$, and that $a$ measures the strength of the 
shock.

\begin{remarque}
\label{f0}
The asymptotic behavior \eqref{asymp} of $j$, $C_1$, and $C_2$ shows that when 
the strength of the shock tends to zero ($a \to 0^+$), the limit of $f(0)$ is 
given by:
\begin{equation*}
f(0) \rightarrow \dfrac{4 \, \gamma^3 \, (\gamma -1)^2}{R^4 \, (\gamma +1)^8} 
\, \dfrac{(c_- +(\gamma-1) \, \frac{e_-}{c_-})^7}{\rho_- \, c_-} >0 \, .
\end{equation*}
\end{remarque}

Since $v_+ <v_-$ for a $1$-shock, we are searching for a solution to \eqref{ode} 
that is defined on all $\R$, and that satisfies:
\begin{equation}
\label{asy-ode}
\lim_{\xi \to -\infty} (\hat{v},w)(\xi) = (a,0) \, ,\quad 
\lim_{\xi \to +\infty} (\hat{v},w)(\xi) = (-a,0) \, .
\end{equation}
To prove Theorem \ref{exis}, we are thus reduced to showing the existence 
of a heteroclinic orbit for \eqref{ode} that connects the stationary solutions 
$(\pm a,0)$. Due to the previous transformation $\hat{v}=v-(v_-+v_+)/2$, if 
$\hat{v}$ takes its values in $[-a,a]$, then $v=\hat{v}+(v_-+v_+)/2$ will 
take its values in the interval $[v_+,v_-]$, and therefore will not vanish.

\begin{remarque}
The system \eqref{ode} is 'singular' at $\hat{v}=0$. Nevertheless, we are 
searching for a smooth solution connecting $(\pm a,0)$, so that $\hat{v}$ 
vanishes in at least one point. Because $w'=\hat{v}''$ should also have a 
limit at this point, a $C^2$ shock profile can exist only if the equation:
\begin{equation*}
w^2 +f(0) \, w +\dfrac{a^2}{2} =0 \, ,
\end{equation*}
has real roots. The corresponding discriminant condition turns out to be 
much less simple than the one found in \cite{kawanishi} for the 'baby model' 
\eqref{Kburgers}. (In particular, $f(0)$ depends on the shock through the 
constants $j$, and $C_1$). This is a first 'nonexplicit' restriction on 
the shock strength to derive the existence of a smooth shock profile.
\end{remarque}

Due to the singular nature of the system \eqref{ode} at $\hat{v}=0$, it is 
more convenient to work on an auxiliary system of ODEs, where the singularity 
has been eliminated (at least formally) thanks to a change of variables. This 
procedure was already used in \cite{kawanishi}. In the next paragraph, we 
shall introduce this auxiliary system, and complete the proof of Theorem 
\ref{exis}.

\subsection{Existence of a heteroclinic orbit}

We begin with a result on an auxiliary system of ODEs, where the singularity 
at $\hat{v}=0$ has been eliminated:

\begin{prop}
\label{lem}
Assume that $\gamma$ satisfies $1<\gamma<(\sqrt{7}+1)/(\sqrt{7}-1)$, and 
consider the following system of ODEs:
\begin{equation}
\label{ode-ref}
\left\{
\begin{array}{rcl}
V' & = & V\, W \, ,\\
W' & = & -W^2-f(V) \, W+\dfrac{(V^2-a^2)}{2} \, .
\end{array}
\right.
\end{equation}
There exists a positive constant $a_0$, that depends only on 
$(\rho_-,u_-,e_-)$, and $\gamma$ such that if the shock strength $a$ 
satisfies $a\in (0,a_0]$, the following properties hold:
\begin{itemize}
\item $f(0)^2-2a^2>0$, and we define $w_0 := \big( 
-f(0)+\sqrt{f(0)^2-2a^2} \big) /2 <0$.

\item There exists a solution $(V_\flat,W_\flat)$ to \eqref{ode-ref} that 
is defined on all $\R$, and that satisfies
\begin{equation*}
\lim_{\eta \rightarrow -\infty} (V_\flat,W_\flat)(\eta) =(a,0) \, ,\quad 
\lim_{\eta \rightarrow +\infty} (V_\flat,W_\flat)(\eta) =(0,w_0) \, .
\end{equation*}
Furthermore, $V_\flat$ is decreasing, and the convergence of $V_\flat$ to 
$0$ as $\eta \rightarrow +\infty$ is exponential.

\item There exists a solution $(V_\sharp,W_\sharp)$ to \eqref{ode-ref} that 
is defined on all $\R$, and that satisfies
\begin{equation*}
\lim_{\eta \rightarrow -\infty} (V_\sharp,W_\sharp)(\eta) =(-a,0) \, ,\quad 
\lim_{\eta \rightarrow +\infty} (V_\sharp,W_\sharp)(\eta) =(0,w_0) \, .
\end{equation*}
Furthermore, $V_\sharp$ is increasing, and the convergence of $V_\sharp$ to 
$0$ as $\eta \rightarrow +\infty$ is exponential.
\end{itemize}
\end{prop}

Assuming that the result of Proposition \ref{lem} holds, the existence of 
a heteroclinic orbit for \eqref{ode} connecting $(\pm a,0)$ can be derived 
by following the analysis of \cite{SchoTad, kawanishi}. We briefly recall the method. 
Using the solution $(V_\flat,W_\flat)$, we introduce the change of variable:
$$
\Xi_\flat (\eta) = -\int_{\eta}^{+\infty} V_\flat (\zeta) \, d\zeta \, .
$$
Since $V_\flat$ tends to $0$ exponentially as $\eta$ tends to $+\infty$, 
$\Xi_\flat$ is well-defined, and it is an increasing $C^\infty$ diffeomorphism 
from $\R$ to $(-\infty,0)$. Then $(\hat{v}_\flat,w_\flat) :=(V_\flat,W_\flat) 
\circ \Xi_\flat^{-1}$ is a $C^\infty$ solution to \eqref{ode} on the interval 
$(-\infty,0)$, that satisfies:
\begin{equation*}
\lim_{\xi \rightarrow -\infty} (\hat{v}_\flat,w_\flat)(\xi) =(a,0) \, ,\quad 
\lim_{\xi \rightarrow 0^-} (\hat{v}_\flat,w_\flat)(\xi) =(0,w_0) \, .
\end{equation*}
Similarly, with the help of the solution $(V_\sharp,W_\sharp)$ we can construct 
a $C^\infty$, decreasing diffeomorphism $\Xi_\sharp$ from $\R$ to $(0,+\infty)$, 
and a $C^\infty$ solution $(\hat{v}_\sharp,w_\sharp)$ to \eqref{ode} on the interval 
$(0,+\infty)$. This solution $(\hat{v}_\sharp,w_\sharp)$ connects $(0,w_0)$ and 
$(-a,0)$, as $\xi$ varies from $0^+$ to $+\infty$. Let us now 'glue' the solutions 
$(\hat{v}_\flat,w_\flat)$, and $(\hat{v}_\sharp,w_\sharp)$, by defining:
\begin{equation}
\label{defsolution}
(\hat{v},w) (\xi) :=\begin{cases}
(\hat{v}_\flat,w_\flat) (\xi) &\text{if $\xi<0$,}\\
(\hat{v}_\sharp,w_\sharp) (\xi) &\text{if $\xi>0$,}
\end{cases}
\end{equation}
and extend the functions $\hat{v}$, and $w$ at $0$ by setting 
$(\hat{v},w)(0)=(0,w_0)$. In this way, $\hat{v}$, and $w$ are continuous 
on $\R$, and $C^\infty$ on $\R \setminus \{ 0\}$. It remains to show that 
$\hat{v} \in C^2(\R)$, that $(\hat{v},w)$ solves \eqref{ode} on $\R$, and 
that $\hat{v}$ takes its values in $(-a,a)$.

Observe first of all that $\hat{v}$ is a decreasing function, because of 
the monotonicity properties of $V_\flat,V_\sharp,\Xi_\flat,\Xi_\sharp$. Using 
the asymptotic behavior of $V_\flat$, $V_\sharp$ at $-\infty$, we get that 
$\hat{v}(\xi) \in (-a,a)$ for all $\xi \in \R$.

Let us now note that the above construction of $(\hat{v},w)$ shows that 
$(\hat{v},w)$ is a solution to \eqref{ode} on $\R \setminus \{ 0\}$. In 
particular, $\hat{v}' (\xi)=w(\xi)$ if $\xi \neq 0$. Moreover, $w$ is 
continuous on $\R$, so $\hat{v} \in C^1(\R)$, and $\hat{v}'(0)=w(0)=w_0$. 
To prove that $\hat{v}\in C^2(\R)$, it is sufficient to show that $w\in 
C^1(\R)$, which is equivalent to showing that $w'$ has a limit at $0$ 
(because we already know that $w$ is $C^\infty$ on $\R \setminus \{ 0\}$). 
To prove that $w'$ has a limit at $0$, we are going to study the asymptotic 
behavior of $(V_\flat,W_\flat)$, and $(V_\sharp,W_\sharp)$ at $+\infty$. 
More precisely, let us denote $U(V,W)$ the vector field associated to the 
ODE \eqref{ode-ref}:
\begin{equation}
\label{defU}
U(V,W)=\begin{pmatrix}
V \, W \\
-W^2- f(V) \, W+\dfrac{(V^2-a^2)}{2} \end{pmatrix} \, ,
\end{equation}
where $f$ is given by \eqref{fun f}. The Jacobian matrix of $U$ at 
$(0,w_0)$ is:
\begin{equation*}
\begin{pmatrix}
w_0 & 0 \\
-f'(0) \, w_0 & -2\, w_0-f(0) \end{pmatrix} = \begin{pmatrix}
\lambda_1^{(0)} & 0 \\
b_0 & \lambda_2^{(0)} \end{pmatrix} \, .
\end{equation*}
For $a$ sufficiently small, one checks that $\lambda_2^{(0)}<\lambda_1^{(0)}<0$ 
(see Proposition \ref{lem} for the definition of $w_0$). The eigenvectors 
corresponding to the eigenvalues $\lambda_1^{(0)}$ and $\lambda_2^{(0)}$ are:
\begin{equation*}
e_1^{(0)}=\begin{pmatrix}
f(0)+3 \, w_0\\
b_0 \end{pmatrix} \, ,\quad e_2^{(0)}=\begin{pmatrix}
0\\ 
1 \end{pmatrix} \, .
\end{equation*}
The standard theory of autonomous ODEs, see e.g. \cite{pontriaguine}, shows that 
there are exactly two solutions to \eqref{ode-ref} that tend to $(0,w_0)$ as 
$\eta$ tends to $+\infty$, and that are tangent to the straight line $(0,w_0) 
+\R \, e_2^{(0)}$. Moreover, all the other solutions to \eqref{ode-ref} that 
tend to $(0,w_0)$ as $\eta$ tends to $+\infty$ are tangent to the straight line 
$(0,w_0) +\R \, e_1^{(0)}$. Now, it is rather simple to see that the two solutions 
to \eqref{ode-ref} that tend to $(0,w_0)$ as $\eta$ tends to $+\infty$, and that 
are tangent to the straight line $(0,w_0) +\R \, e_2^{(0)}$, satisfy 
$V \equiv 0$, and:
\begin{equation*}
W'=-W^2-f(0) \, W -\dfrac{a^2}{2} \, .
\end{equation*}
Because the solutions $(V_\flat,W_\flat)$, and $(V_\sharp,W_\sharp)$ given by 
Proposition \ref{lem} cannot satisfy $V_\flat \equiv 0$, and $V_\sharp \equiv 0$, 
we can conclude that the solutions $(V_\flat,W_\flat)$, and $(V_\sharp,W_\sharp)$ 
are tangent to $(0,w_0) +\R \, e_1^{(0)}$ as $\eta$ tends to $+\infty$. In 
particular, this yields:
\begin{equation}
\label{asympderiv}
\lim_{\eta \rightarrow +\infty} \dfrac{W_\flat' (\eta)}{V_\flat' (\eta)} 
=\lim_{\eta \rightarrow +\infty} \dfrac{W_\sharp' (\eta)}{V_\sharp' (\eta)} 
=\dfrac{-f'(0) \, w_0}{f(0)+3 \, w_0} \, .
\end{equation}
(A quick verification shows that $f(0)+3 \, w_0>0$ for small enough $a$.) 
From the construction of the solutions $(\hat{v}_\flat,w_\flat)$, and 
$(\hat{v}_\sharp,w_\sharp)$, we get:
\begin{equation*}
\lim_{\xi \rightarrow 0^-} w_\flat' (\xi) 
=\lim_{\xi \rightarrow 0^+} w_\sharp' (\xi) 
=\dfrac{-f'(0) \, w_0^2}{f(0)+3 \, w_0} \, .
\end{equation*}
As a consequence, when $a$ is small enough, $w \in C^1(\R)$, and therefore 
$\hat{v} \in C^2(\R)$. Moreover, $(\hat{v},w)$ solves \eqref{ode} on 
$\R \setminus \{ 0\}$, so by continuity, it solves \eqref{ode} on $\R$. 
This completes the proof of Theorem \ref{exis}, provided that the result 
of Proposition \ref{lem} holds.

\subsection{Proof of Proposition \ref{lem}}

\begin{figure}
\centering
\psfrag{A}{$a$}
\psfrag{B}{$-a$}
\psfrag{C}{$V$}
\psfrag{D}{$W$}
\psfrag{E}{$\overline{V}$}
\psfrag{P1}{$P_1$}
\psfrag{P2}{$P_2$}
\psfrag{w0}{$w_0$}
\includegraphics{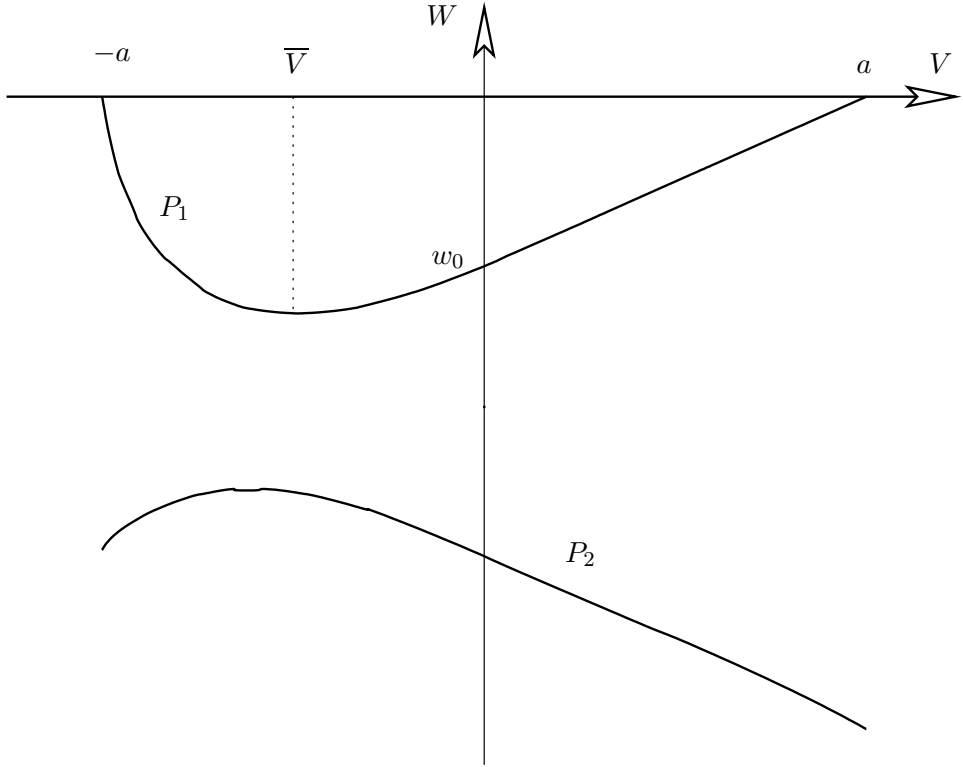}
\caption{The set $P =P_1 \cup P_2$.}
\label{K}
\end{figure}

In this paragraph, we prove Proposition \ref{lem}, which will complete the 
proof of Theorem \ref{exis}. At first, we define the set:
$$
P = \Big\{ (V,W) | V \in [-a,a], W^2+f(V) \, W-\dfrac{V^2-a^2}{2}=0 \Big\} \, ,
$$
so that the points $(\pm a,0)$ belong to $P$. The following Lemma gives a 
description of $P$ for $a>0$ small enough. We refer to figure \ref{K} for 
a schematic picture.

\begin{lemm}
\label{geo}
Assume that $1 < \gamma < (\sqrt{7}+1)/(\sqrt{7}-1)$. Then there exists 
a constant $a_0>0$, that only depends on $(\rho_-,u_-,e_-)$ and $\gamma$ 
such that if the shock strength $a$ satisfies $a\in (0,a_0]$, we have the 
following results:
\begin{itemize}
\item For all $V \in [-a,a]$, $f(V)^2+2\, (V^2-a^2)>0$. We can 
thus define
\begin{align*}
\forall \, V \in [-a,a] \, ,\quad 
\W_1(V) &= \dfrac{-f(V)+\sqrt{f(V)^2 +2\, (V^2-a^2)}}{2} \, ,\\
\W_2(V) &= \dfrac{-f(V)-\sqrt{f(V)^2 +2\, (V^2-a^2)}}{2} \, .
\end{align*}

\item $P = P_1\cup P_2$, where $P_1$ and $P_2$ are two curves defined by
$$
P_1 =\big\{ (V,\W_1(V)) | V \in [-a,a] \big\} \, , \quad 
P_2 =\big\{ (V,\W_2(V)) | V \in [-a,a] \big\} \, ,
$$
so that the points $(\pm a,0)$, and $(0,w_0)$ belong to $P_1$.

\item There exists a unique point $\overline{V} \in \, (-a,0)$ such that 
$\W_1$ is increasing on the interval $[\overline{V},a]$, and $\W_1$ is 
decreasing on the interval $[-a,\overline{V}]$.

\item For all $V \in [\overline{V},0]$, one has 
$\W_2(V) < \W_1(\overline{V})$.
\end{itemize}
\end{lemm}

\begin{proof}
Let us first define a function $\Delta$ by setting:
\begin{equation*}
\Delta (V) := f(V)^2 +2 \, (V^2-a^2) \, .
\end{equation*}
Using \eqref{asymp}, for $a$ small enough, we have:
\begin{equation*}
\dfrac{C_1}{\gamma +1} -a \ge \kappa >0 \, ,\quad 
\dfrac{\gamma \, C_1}{\gamma +1} -a \ge \kappa >0 \, ,\quad 
\dfrac{(\gamma -1)\, C_1}{\gamma +1} -2a \ge \kappa >0 \, ,
\end{equation*}
where $\kappa$ is a positive constant that only depends on $(\rho_-,u_-,e_-)$ 
and $\gamma$. Moreover, \eqref{asymp} also shows that $j \ge \kappa$ for 
$a \in (0,a_0]$, up to restricting $\kappa$. Consequently, there exist 
$a_0>0$, and $\kappa>0$ such that for $a \in (0,a_0]$, we have $f(V) \ge \kappa$, 
and $\Delta (V) \ge \kappa$ for all $V\in [-a,a]$. This directly shows that 
the set $P$ is the union of the two curves $P_1$, and $P_2$. It is rather 
clear from the definition of $P_1$ that $(\pm a,0)$, and $(0,w_0)$ belong 
to $P_1$ (recall that $w_0$ is defined in Proposition \ref{lem}). Observe 
also that $\W_2(V) <\W_1 (V) \le 0$ for $V \in [-a,a]$, and $\W_1(V)<0$ if 
$V \in (-a,a)$.

The functions $\W_1$, and $\W_2$ are $C^\infty$ on $[-a,a]$. Moreover, we 
compute the relation:
\begin{equation}
\label{derivW1}
\forall \, V \in [-a,a] \, ,\quad 
\sqrt{\Delta (V)} \, \W_1'(V) =V-\W_1(V) \, f'(V) \, ,
\end{equation}
and from \eqref{fun f}, we also compute
\begin{multline}
\label{derivf}
f'(V)=\dfrac{14 \, (\gamma-1)}{j \, R^4 \, (\gamma+1)} \\
\left( \dfrac{C_1}{\gamma+1}-\hat{v} \right)^2 \, 
\left( \hat{v}+\dfrac{\gamma\, C_1}{\gamma+1} \right)^2 \, 
\left( 2\hat{v}+\dfrac{(\gamma-1) \, C_1}{\gamma+1} +\dfrac{C_1}{\sqrt{7}} \right) 
\left( -2\hat{v}-\dfrac{(\gamma-1) \, C_1}{\gamma+1} +\dfrac{C_1}{\sqrt{7}} \right) \, .
\end{multline}
As we have done for $f$, and $\Delta$, a careful analysis shows that for 
$1<\gamma<(\sqrt{7}+1)/(\sqrt{7}-1)$, and for $a$ small enough, one has 
$f'(V) \ge \kappa>0$ for all $V \in [-a,a]$, because each term in the product 
\eqref{derivf} is positive. Using this information in \eqref{derivW1}, we can 
already conclude that $\W_1$ is increasing on the interval $[0,a]$ (see figure 
\ref{K}). Moreover, the relation \eqref{derivW1} also shows that $\W'_1(0)>0$, 
and $\W'_1(-a)<0$. Consequently, there exists some $\overline{V} \in (-a,0)$ 
such that $\W_1'(\overline{V})=0$. Let us prove that $\overline{V}$ is the only 
zero of $\W_1'$. We claim that it is sufficient to show the following property:
\begin{equation}
\label{property}
\W_1' (V) =0 \Longrightarrow \W_1'' (V) >0 \, .
\end{equation}
Indeed, if the property \eqref{property} holds true, then any point where $\W_1'$ 
vanishes is a local strict minimum. If there existed two such local strict minima 
$-a<\overline{V}_1<\overline{V}_2<a$, then $\W_1$ would admit a local maximum 
$\overline{V}_3 \in (\overline{V}_1,\overline{V}_2)$, which is obviously 
impossible. Therefore let us prove that the property \eqref{property} holds 
true.

Differentiating \eqref{derivW1} with respect to $V$, we obtain that if 
$\W_1'(\overline{V})=0$, then
$$
\sqrt{\Delta (\overline{V})} \, \W_1''(\overline{V}) = 1 -f''(\overline{V}) 
\, \W_1(\overline{V}) \, .
$$
Observing that
\begin{equation*}
|f''(\overline{V}) \, \W_1(\overline{V})| \le C \, |\W_1(\overline{V})| 
\le C \, \dfrac{a^2 -\overline{V}^2}{f(\overline{V})} \le 
\dfrac{C \, a^2}{\kappa} \, ,
\end{equation*}
for suitable positive constants $C$, and $\kappa$ (that are independent 
of $a \in (0,a_0]$), we can conclude that $\W_1''(\overline{V})>0$, provided 
that $a$ is small enough. This completes the proof that $\W_1'$ has a unique 
zero $\overline{V} \in (-a,0)$, and therefore $\W_1$ is decreasing on 
$[-a,\overline{V}]$, and is increasing on $[\overline{V},a]$.

For the last point of the lemma, we use the relation:
$$
\W_1'(V) +\W_2'(V)=-f'(V)<0 \, .
$$
Because $\W_1'(V)\ge 0$ for $V \in [\overline{V},0]$, we get $\W_2'(V)<0$ 
for $V \in [\overline{V},0]$. Thus for $V \in [\overline{V},0]$, we have 
$\W_2(V) \le \W_2 (\overline{V}) <\W_1 (\overline{V})$, and the proof of 
the Lemma is complete.
\end{proof}

Using Lemma \ref{geo}, we are going to prove Proposition \ref{lem}. 
The analysis follows \cite{gilbarg}.

\begin{figure}
\centering
\psfrag{$W$}{$W$}
\psfrag{$V$}{$V$}
\psfrag{$w_0$}{$w_0$}
\psfrag{$a$}{$a$}
\psfrag{$K_1$}{$K_1$}
\psfrag{A}{$(a,0)+{\mathbb R} \, r_2$}
\includegraphics{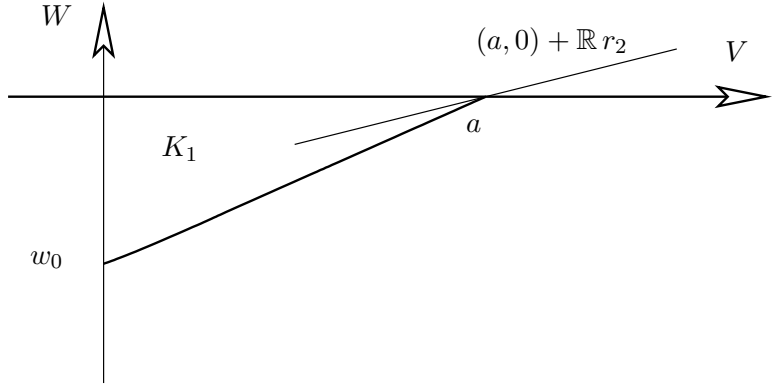}
\caption{The compact set $K_1$.}
\label{compactK1}
\end{figure}

As we have already seen in the preceeding paragraph, the point $(w_0,0)$ is 
a stable node of \eqref{ode-ref}. We now study the nature of the equilibrium 
points $(\pm a,0)$. Recall that the vector field associated to \eqref{ode-ref} 
is denoted $U$, see \eqref{defU}. The Jacobian matrix of $U$ at $(a,0)$ is:
\begin{equation*}
\begin{pmatrix}
0 & a \\
a & f(a) \end{pmatrix} \, ,
\end{equation*}
so it has exactly one negative eigenvalue $\mu_1$, and one positive eigenvalue 
$\mu_2$ (the equilibrium point $(a,0)$ is a saddle point):
\begin{equation*}
\mu_1 =\dfrac{-f(a)-\sqrt{f(a)^2+4 \, a^2}}{2} \, ,\quad 
\mu_2 =\dfrac{-f(a)+\sqrt{f(a)^2+4 \, a^2}}{2} \, .
\end{equation*}
An eigenvector associated to $\mu_2$, and is $r_2=(a,\mu_2)$. Moreover, 
using the relation \eqref{derivW1}, we can check that for $a$ small enough, 
the following inequality holds:
\begin{equation}
\label{ineg1}
0 <\dfrac{\mu_2}{a} < \dfrac{a}{f(a)} =\W_1'(a) \, ,
\end{equation}
where the function $\W_1$ is defined in Lemma \ref{geo}. Let us now define a 
compact set $K_1$ by:
\begin{equation*}
K_1 := \Big\{ (V,W) \in [0,a] \times \R \, | \, \W_1 (V) \le W \le 0 \Big\} \, ,
\end{equation*}
Then the inequalities \eqref{ineg1} show that for $s<0$ small enough, the 
point $(a,0) + s \, r_2$ belongs to the interior of $K_1$. We refer to 
figure \ref{compactK1} for a detailed picture of the situation.

From the standard theory of autonomous ODEs, see e.g. \cite{pontriaguine}, we 
know that there exists a maximal solution $(V_\flat,W_\flat)$ to \eqref{ode-ref} 
that tends to the saddle point $(a,0)$ as $\eta$ tends to $-\infty$, and that 
is tangent to the half-straight line $(a,0) +\R^- \, r_2$. This solution is 
defined on an open interval $(-\infty,\eta_*)$ (with possibly $\eta_*=+\infty$). 
For large negative $\eta$, the preceeding analysis shows that 
$(V_\flat,W_\flat) (\eta)$ belongs to the interior of $K_1$. Moreover, 
$(V_\flat,W_\flat)$ cannot reach the boundary of $K_1$. Indeed $V_\flat$ cannot 
identically vanish so $(V_\flat,W_\flat)(\eta) \not \in \partial K_1 \cap 
\{ V=0\}$. Similarly, we have $(V_\flat,W_\flat)(\eta) \neq (a,0)$. Eventually, 
on the set:
\begin{equation*}
\Big\{ (V,0) \, | \, V \in (0,a) \Big\} \cup 
\Big\{ (V,\W_1(V)) \, | \, V \in \, (0,a) \Big\} \, ,
\end{equation*}
the vector field $U$ is not zero, and is directed towards the interior of 
$K_1$. Therefore the solution $(V_\flat,W_\flat)$ cannot reach $\partial K_1$, 
so it takes its values in the compact set $K_1$. The maximal solution 
$(V_\flat,W_\flat)$ is thus defined on $\R$. It cannot reach the boundary 
of $K_1$, so $W_\flat$ takes negative values, which means that $V_\flat$ is 
decreasing (because $V_\flat'=V_\flat \, W_\flat$). Because $(V_\flat,W_\flat)$ 
takes values in the interior of $K_1$, the function $W_\flat$ is also decreasing. 
This shows that $(V_\flat,W_\flat)(\eta)$ has a limit as $\eta$ tends to $+\infty$, 
and this limit is necessarily  be a stationary solution of \eqref{ode-ref}. The 
only possibility is that $(V_\flat,W_\flat)(\eta)$ tends to $(0,w_0)$ as $\eta$ 
tends to $+\infty$. The convergence is necessarily exponential, because the 
Jacobian matrix of $U$ at $(0,w_0)$ has two negative eigenvalues, see e.g. 
\cite{pontriaguine}.

\begin{figure}
\centering
\psfrag{$W$}{$W$}
\psfrag{$V$}{$V$}
\psfrag{$K_2$}{$K_2$}
\psfrag{$-a$}{$-a$}
\psfrag{A}{$\overline{V}$}
\psfrag{B}{$W_1(\overline{V})$}
\psfrag{C}{$(-a,0)+{\mathbb R} \, R_2$}
\includegraphics{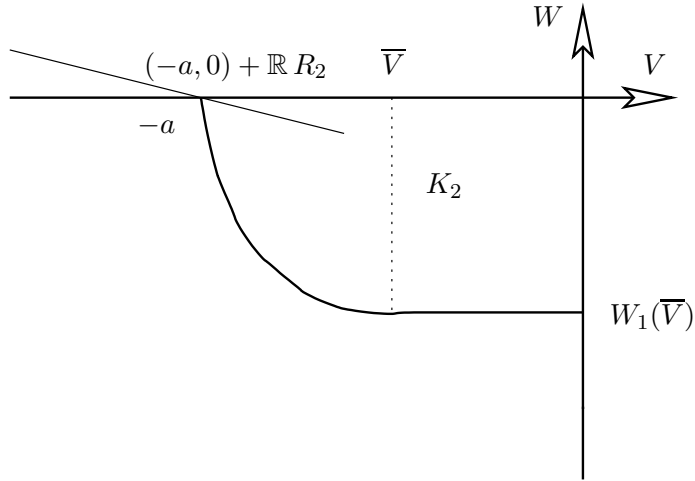}
\caption{The compact set $K_2$.}
\label{compactK2}
\end{figure}

To construct the other solution $(V_\sharp,W_\sharp)$, we argue similarly by defining 
a compact set $K_2$:
\begin{equation*}
K_2 := \Big\{ (V,W) \in [-a,\overline{V}] \times \R | \W_1 (V) \le W \le 0 \Big\} 
\cup \Big\{ 
(V,W) \in [\overline{V},0] \times \R | \W_1 (\overline{V}) \le W \le 0 \Big\} \, ,
\end{equation*}
see figure \ref{compactK2}. The Jacobian matrix of $U$ at $(-a,0)$ has one 
negative eigenvalue $\nu_1$, and one positive eigenvalue $\nu_2$, with:
\begin{equation*}
\nu_2 = \dfrac{-f(-a)+\sqrt{f(-a)^2+4 \, a^2}}{2} \, .
\end{equation*} 
An eigenvector associated to the eigenvalue $\nu_2$ is $R_2=(-a,\nu_2)$. 
As was done earlier, we check that the inequalities:
\begin{equation*}
\W_1'(-a)=\dfrac{-a}{f(-a)} < \dfrac{\nu_2}{-a} <0 \, ,
\end{equation*}
hold true. Therefore, one can reproduce the above analysis, and show that there 
exists a solution $(V_\sharp,W_\sharp)$ to \eqref{ode-ref} that takes its values 
in $K_2$ (and is thus defined on $\R$), and that tends to $(-a,0)$ at $-\infty$. 
Moreover, $(V_\sharp,W_\sharp)$ can not reach the boundary of $K_2$, so $V_\sharp$ 
is increasing. It only remains to study the monotonicity of $W_\sharp$. This is 
slightly more complicated than for $W_\flat$. Observe that $K_2$ is the union of 
the sets:
\begin{align*}
K_2^1 &:=\Big\{ (V,W) \in [-a,0] \times \R \, | \, \W_1 (V) \le W \le 0 \Big\} \, ,\\
K_2^2 &:=\Big\{ (V,W) \in [\overline{V},0] \times \R \, | \, 
\W_1 (\overline{V}) \le W \le \W_1(V) \Big\} \, .
\end{align*}
When $(V_\sharp,W_\sharp)$ takes its values in the interior of $K_2^1$, the function 
$W_\sharp$ is decreasing (this is the case for large negative $\eta$). At the opposite, 
when $(V_\sharp,W_\sharp)$ takes its values in the interior of $K_2^2$, the function 
$W_\sharp$ is increasing, because thanks to Lemma \ref{geo}, we have:
\begin{align*}
W_\sharp'(\eta) 
&=-W_\sharp(\eta)^2 -f(V_\sharp(\eta)) \, W_\sharp(\eta) 
+\dfrac{V_\sharp(\eta)^2 -a^2}{2} \\
&=\big( 
\W_1 (V_\sharp(\eta)) -W_\sharp (\eta) \big) \, 
\big( W_\sharp(\eta) -\W_2 (V_\sharp(\eta)) \big) \\
&\ge \big( \W_1 (V_\sharp(\eta)) -W_\sharp (\eta) \big) \, 
\big( \W_1(\overline{V}) -\W_2 (V_\sharp(\eta)) \big) >0 \, .
\end{align*}
Moreover, if $(V_\sharp,W_\sharp) (\eta_0)$ belongs to the interior of $K_2^2$ for 
some $\eta_0 \in \R$, then $(V_\sharp,W_\sharp) (\eta)$ belongs to the interior of 
$K_2^2$ for all $\eta \ge \eta_0$ (because it cannot reach the boundary of $K_2^2$ 
for $\eta \ge \eta_0$). Summing up, either $(V_\sharp,W_\sharp)(\eta)$ belongs to 
$K_2^1$ for all $\eta$, and $W_\sharp$ is monotonic on $\R$, either 
$(V_\sharp,W_\sharp)(\eta)$ belongs to $K_2^2$ for all $\eta$ greater than some 
$\eta_0$, and $W_\sharp$ is monotonic on $[\eta_0,+\infty)$. In any case, the function 
$W_\sharp$ is monotonic on a neighborhood of $+\infty$, and thus has a limit at 
$+\infty$. This shows that $(V_\sharp,W_\sharp)(\eta)$ tends to $(0,w_0)$ as $\eta$ 
tends to $+\infty$, and the convergence is exponential. As a matter of fact, we 
have seen in the preceeding paragraph that $(V_\sharp,W_\sharp)$ is tangent to the 
straight line $(0,w_0)+\R \, e_1^{(0)}$ as $\eta$ tends to $+\infty$, so one can 
check that $(V_\sharp,W_\sharp) (\eta)$ belongs to the interior of $K_2^2$ for 
large positive $\eta$. This means that $W_\sharp$ is decreasing on some interval 
$(-\infty,\eta_0)$, and increasing on $[\eta_0,+\infty)$. The proof of 
Proposition \ref{lem} is now complete.

\section{Additional regularity of shock profiles}
\label{moresmooth}

As should be clear from the preceeding section, the key point in the construction 
of a shock profile is Proposition \ref{lem} that gives the existence of two 
heteroclinic orbits for the system \eqref{ode-ref}. To prove Theorem \ref{smooth}, 
we are going to study the behavior of the derivatives of $(V_\flat,W_\flat)$, and 
$(V_\sharp,W_\sharp)$ near $+\infty$. The proof of Theorem \ref{smooth} follows 
from an induction argument. To make the arguments clear, we deal with the first 
case separately. In all what follows, $(V_\flat,W_\flat)$, and $(V_\sharp,W_\sharp)$ 
are the solutions to \eqref{ode-ref} that are defined in Proposition \ref{lem}, 
and $(\hat{v},w)$ denotes the solution to \eqref{ode} that is defined by 
\eqref{defsolution}. We have the following:

\begin{prop}
\label{n=3}
Under the assumptions of Proposition \ref{lem}, there exists a positive constant 
$a_1 \le a_0$ (that depends on $(\rho_-,u_-,e_-)$, and $\gamma$), such that for 
all $a \in (0,a_1]$, one has $w \in C^2 (\R)$, $\hat{v} \in C^3 (\R)$, and:
\begin{equation*}
w(\xi)=w_0 +w_1 \, \hat{v} (\xi) +w_2 \, \hat{v} (\xi)^2 +o(\hat{v} (\xi)^2) \, ,\quad 
\text{\rm as $\xi \rightarrow 0$,}
\end{equation*}
for some suitable constants $w_1,w_2$ ($w_0$ has already been defined in Proposition 
\ref{lem}).
\end{prop}

\begin{proof}
Recall that $V_b$, and $V_\sharp$ do not vanish on $\R$, so we can introduce 
some $C^\infty$ functions $W_{\flat,1}$, and $W_{\sharp,1}$ that are defined by:
\begin{equation*}
W_\flat = w_0 + V_\flat \, W_{\flat,1} \, ,\quad 
W_\sharp = w_0 + V_\sharp \, W_{\sharp,1} \, .
\end{equation*}
Substituting in \eqref{ode-ref} shows that $(V_\flat,W_{\flat,1})$, and 
$(V_\sharp,W_{\sharp,1})$ are solutions to the system:
\begin{equation}
\label{w_1}
\begin{cases}
V' = V (w_0+V \, W_1) \, , &\\
W_1' = -w_0 \, \dfrac{f(V)-f(0)}{V} -2 \, V \, W_1^2 -3 \, w_0 \, W_1 
-f(V) \, W_1 +\dfrac{V}{2} \, . &
\end{cases}
\end{equation}
Moreover, we already know from Proposition \ref{lem}, and \eqref{asympderiv} 
that:
\begin{equation*}
\lim_{\eta \rightarrow +\infty} (V_\flat,W_{\flat,1}) (\eta) 
=\lim_{\eta \rightarrow +\infty} (V_\sharp,W_{\sharp,1}) (\eta) 
=\Big( 0,\dfrac{-f'(0) \, w_0}{f(0)+3 \, w_0} \Big) \, .
\end{equation*}

We denote $U_1 (V,W_1)$ the vector field associated with \eqref{w_1}:
\begin{equation*}
U_1 (V,W_1) := \begin{pmatrix} 
V (w_0+V \, W_1) \\
-w_0 \, \dfrac{f(V)-f(0)}{V} -2 \, V \, W_1^2 -3 \, w_0 \, W_1 
-f(V) \, W_1 +\dfrac{V}{2} \end{pmatrix} \, .
\end{equation*}
Recall that $f$ is a polynomial function of degree $7$, see \eqref{fun f}, 
thus $F(V) :=(f(V)-f(0))/V$ is a polynomial function of degree $6$, and we have 
$F(0)=f'(0)$, $F'(0)=f''(0)/2$. Obviously the system of ODEs \eqref{w_1} admits 
the equilibrium point $(0,w_1)$, where:
\begin{equation*}
w_1 :=\dfrac{-F(0) \, w_0}{f(0)+3\, w_0} =\dfrac{-f'(0) \, w_0}{f(0)+3\, w_0} \, .
\end{equation*}
We are now going to study the nature of the equilibrium point $(0,w_1)$, and show 
that for $a$ small enough, this equilibrium point is a stable node for \eqref{w_1}. 
Then we shall show that $w \in C^2 (\R)$, and $\hat{v} \in C^3 (\R)$. In the end, 
we shall derive the asymptotic expansion near $\xi=0$.

\underline{Step 1}: the Jacobian matrix of $U_1$ at $(0,w_1)$ is:
\begin{equation*}
\begin{pmatrix}
w_0 & 0\\
\dfrac{1}{2} -2\, w_1^{2} -f'(0)\, w_1 -\dfrac{f''(0)}{2} \, w_0 & -f(0)-3\, w_0
\end{pmatrix} =\begin{pmatrix}
\lambda_1^{(1)} & 0 \\
b_1 & \lambda_2^{(1)} \end{pmatrix} \, .
\end{equation*}
Using Remak \ref{f0}, we can conclude that for sufficiently small $a$, that 
is $a \in (0,a_1]$ for some positive number $a_1$ less than $a_0$, one has 
$\lambda_2^{(1)} <\lambda_1^{(1)}<0$, that is, $f(0)+4\, w_0>0$. Moreover, 
the eigenvectors corresponding to the eigenvalues $\lambda_1^{(1)}$, and 
$\lambda_2^{(1)}$ are:
\begin{equation*}
e_1^{(1)} =\begin{pmatrix}
f(0) +4\, w_0\\ 
b_1 \end{pmatrix} \, ,\quad
e_2^{(1)} =\begin{pmatrix}
0\\ 
1 \end{pmatrix} \, .
\end{equation*}
Consequently, $(0,w_1)$ is a stable node of \eqref{w_1}, and there are exactly 
two solutions to \eqref{w_1} that tend to $(0,w_1)$ as $\eta$ tends to $+\infty$, 
and that are tangent to the straight line $(0,w_1) +\R \, e_2^{(1)}$. All the other 
solutions to \eqref{w_1} that tend to $(0,w_1)$ as $\eta$ tends to $+\infty$ are 
tangent to the straight line $(0,w_1) +\R \, e_1^{(1)}$. As in the preceeding 
section, we can thus conclude that:
\begin{equation}
\label{eq:0}
\lim_{\eta \rightarrow +\infty} \dfrac{W_{\flat,1}'(\eta)}{V_\flat'(\eta)} 
=\lim_{\eta \rightarrow +\infty} \dfrac{W_{\sharp,1}'(\eta)}{V_\sharp'(\eta)} 
=\dfrac{b_1}{f(0)+4\, w_0} \, .
\end{equation}

\underline{Step 2}: if we let $g_1$ denote the second coordinate of 
the vector field $U_1$, we have $W_{\flat,1}'=g_1(V_\flat,W_{\flat,1})$, and 
$W_{\sharp,1}'=g_1(V_\sharp,W_{\sharp,1})$. Differentiating once with respect 
to $\eta$, and using \eqref{eq:0}, we end up with:
\begin{equation}
\label{eq:1}
\lim_{\eta \rightarrow +\infty} \dfrac{W_{\flat,1}''(\eta)}{V_\flat'(\eta)} 
=\lim_{\eta \rightarrow +\infty} \dfrac{W_{\sharp,1}''(\eta)}{V_\sharp'(\eta)} 
= \ell_1 \, ,
\end{equation}
where the real number $\ell_1$ can be explicitely computed (but its exact expression 
is of no use). Following the analysis of the preceeding section, we define some 
functions $w_{\flat,1} := W_{\flat,1} \circ \Xi_\flat^{-1}$, and $w_{\sharp,1} 
:= W_{\sharp,1} \circ \Xi_\sharp^{-1}$. First of all, \eqref{eq:0} yields:
\begin{equation}
\label{eq:2}
\lim_{\xi \rightarrow 0^-} w_{\flat,1}' (\xi) 
=\lim_{\xi \rightarrow 0^+} w_{\sharp,1}' (\xi) 
=\dfrac{b_1 \, w_0}{f(0)+4\, w_0} \, .
\end{equation}
Observe now that we have the relations:
\begin{equation*}
\hat{v}_\flat \, w_{\flat,1}' =W_{\flat,1}' \circ \Xi_\flat^{-1} \, ,\quad 
\hat{v}_\sharp \, w_{\sharp,1}' =W_{\sharp,1}' \circ \Xi_\sharp^{-1} \, ,
\end{equation*}
and combining with \eqref{eq:1}, we get:
\begin{equation}
\label{eq:3}
\begin{split}
\lim_{\xi \rightarrow 0^-} (\hat{v}_\flat \, w_{\flat,1}')' (\xi) 
=\lim_{\eta \rightarrow +\infty} \dfrac{W_{\flat,1}''(\eta)}{V_\flat (\eta)} 
=\lim_{\eta \rightarrow +\infty} 
\dfrac{W_{\flat,1}''(\eta) \, W_\flat (\eta)}{V_\flat' (\eta)} =\ell_1 \, w_0 \, ,\\
\lim_{\xi \rightarrow 0^+} (\hat{v}_\sharp \, w_{\sharp,1}')' (\xi) 
=\lim_{\eta \rightarrow +\infty} \dfrac{W_{\sharp,1}''(\eta)}{V_\sharp (\eta)} 
=\lim_{\eta \rightarrow +\infty} 
\dfrac{W_{\sharp,1}''(\eta) \, W_\sharp (\eta)}{V_\sharp' (\eta)} =\ell_1 \, w_0 \, .
\end{split}
\end{equation}
Differentiating twice the relations $w_\flat =w_0 +\hat{v}_\flat \, w_{\flat,1}$, 
and $w_\sharp =w_0 +\hat{v}_\sharp \, w_{\sharp,1}$, we obtain:
\begin{gather*}
w_\flat'' =\hat{v}_\flat'' \, w_{\flat,1} +\hat{v}_\flat' \, w_{\flat,1}' 
+(\hat{v}_\flat \, w_{\flat,1}')' =w_\flat' \, w_{\flat,1} 
+w_\flat \, w_{\flat,1}' +(\hat{v}_\flat \, w_{\flat,1}')' \, ,\\
w_\sharp'' =w_\sharp' \, w_{\sharp,1} 
+w_\sharp \, w_{\sharp,1}' +(\hat{v}_\sharp \, w_{\sharp,1}')' \, .
\end{gather*}
Using \eqref{eq:2}, and \eqref{eq:3}, we get $w_\flat'' (0^-)=w_\sharp'' (0^+)$. 
Using the definition \eqref{defsolution}, this shows that $w \in C^2 (\R)$, and 
using $\hat{v}'=w$, we obtain $\hat{v} \in C^3 (\R)$.

\underline{Step 3}: note that we have the following expansions near $\xi=0$:
\begin{align*}
&w(\xi) = w(0) +w'(0) \, \xi +\dfrac{w''(0)}{2} \, \xi^2 +o(\xi^2) \, ,\\
&\hat{v}(\xi) = w(0) \, \xi +\dfrac{w'(0)}{2} \, \xi^2 +o(\xi^2) \, ,
\end{align*}
with $w(0)=w_0 <0$. We can thus combine these expansions, and derive:
\begin{equation*}
w(\xi) = w_0 +\alpha \, \hat{v}(\xi) +\beta \, \hat{v}(\xi)^2 
+o(\hat{v}(\xi)^2) \, ,
\end{equation*}
for some appropriate real numbers $\alpha$, and $\beta$, that we are going 
to determine. From the relation $w_\flat (\xi)=w_0+\hat{v}_\flat(\xi) \, 
w_{\flat,1} (\xi)$, and using that $w_{\flat,1} (\xi)$ tends to $w_1$ as 
$\xi$ tends to $0^-$, we first get $\alpha=w_1$. Then from \eqref{eq:2}, 
and from the relation $\hat{v}_\flat'(0^-)=w_0$, we can obtain:
\begin{equation*}
w_{\flat,1} (\xi)=w_1 +\dfrac{b_1}{f(0)+4\, w_0} \, \hat{v}_\flat (\xi) 
+o(\hat{v}_\flat (\xi)) \, ,\quad \text{as $\xi \rightarrow 0^-$.}
\end{equation*}
We thus obtain $\beta=b_1/(f(0)+4\, w_0)$, which yields:
\begin{equation*}
w(\xi) = w_0 +w_1 \, \hat{v}(\xi) +w_2 \, \hat{v}(\xi)^2 +o(v(\xi)^2) \, ,
\end{equation*}
where $w_2 :=b_1/(f(0)+4\, w_0)$. This latter expansion will be generalized 
to any order in what follows.
\end{proof}

We now turn to the proof of Theorem \ref{smooth}. More precisely, we are going 
to prove the following result, that is a refined version of Theorem \ref{smooth}:

\begin{thh}
\label{smooth2}
Let the assumptions of Proposition \ref{lem} be satisfied. Then there exists a 
nonincreasing sequence of positive numbers $(a_n)_{n \in \N}$ such that, for 
all integer $n$, if $a \in (0,a_n]$, then $w \in C^{n+1} (\R)$, and $\hat{v} 
\in C^{n+2} (\R)$. Moreover, $w$ admits the following asymptotic expansion 
near $\xi=0$:
\begin{equation}
\label{expansion}
w(\xi) =w_0+w_1 \, \hat{v}(\xi) +\cdots +w_{n+1} \, \hat{v} (\xi)^{n+1} 
+o(\hat{v} (\xi)^{n+1}) \, ,
\end{equation}
where the real numbers $w_0,\dots,w_{n+1}$ are defined by:
\begin{equation*}
\begin{cases}
w_0 =\dfrac{-f(0)+\sqrt{f(0)^2 -2\, a^2}}{2} \, ,& \\
w_k =\dfrac{b_{k-1}}{f(0) +(k+2)\, w_0} \, ,&\text{\rm for $k=1,\dots,n+1$,}
\end{cases}
\end{equation*}
and the real numbers $b_0,\dots,b_n$ are given by:
\begin{equation*}
\begin{cases}
b_0 = -f'(0) \, w_0 \, ,& \\
b_1 =\dfrac{1}{2}-2 \, w_1^2 -f'(0) \, w_1 -\dfrac{f''(0)}{2} \, w_0 \, ,& \\
b_k =-\displaystyle \sum_{i=1}^{k+1} \dfrac{f^{(i)}(0)}{i!} \, w_{k+1-i} 
-\displaystyle \sum_{i=1}^k (i+1) \, w_i \, w_{k+1-i} \, ,& 
\text{\rm for $k=2,\dots,n$.}
\end{cases}
\end{equation*}
\end{thh}

\begin{proof}
The case $n=0$ has been proved in the preceeding section, while the case $n=1$ 
is proved in Proposition \ref{n=3}. (The reader can check that the definition of 
$w_0$, $w_1$, $w_2$, $b_0$, and $b_1$ coincide with our previous notations.) We 
prove the general case by using an induction with respect to $n$, and we thus 
assume that the result of Theorem \ref{smooth2} holds up to the order $n \ge 1$. 
We are going to construct $a_{n+1}$ so that the conclusion of Theorem \ref{smooth2} 
holds for $a \in (0,a_{n+1}]$. In particular, the real numbers $w_0,\dots,w_{n+1}$, 
and $b_0,\dots,b_n$ are given as in Theorem \ref{smooth2}, and we can already define 
the real number $b_{n+1}$ by the formula:
\begin{equation*}
b_{n+1} :=-\displaystyle \sum_{i=1}^{n+2} \dfrac{f^{(i)}(0)}{i!} \, w_{n+2-i} 
-\displaystyle \sum_{i=1}^{n+1} (i+1) \, w_i \, w_{n+2-i} \, .
\end{equation*}
(Observe indeed that this definition only involves $w_0,\dots,w_{n+1}$, and not 
$w_{n+2}$.)

\underline{Step 1}: because $V_\flat$, and $V_\sharp$ do not vanish, we can 
introduce some functions $W_{\flat,n+1}$, and $W_{\sharp,n+1}$ by the relations:
\begin{equation*}
W_\flat =w_0 +w_1 \, V_\flat +\cdots +w_n \, V_\flat^n 
+W_{\flat,n+1} \, V_\flat^{n+1} \, ,\quad 
W_\sharp =w_0 +w_1 \, V_\sharp +\cdots +w_n \, V_\sharp^n 
+W_{\sharp,n+1} \, V_\sharp^{n+1} \, .
\end{equation*}
Thanks to Taylor's formula, we can write the polynomial function $f$ as:
\begin{equation*}
f(V)=f(0)+f'(0) \, V+\dfrac{f''(0)}{2} \, V^2 +\cdots 
+\dfrac{f^{(n)}(0)}{n!} \, V^n +V^{n+1} \, F_{n+1} (V) \, ,
\end{equation*}
where $F_{n+1}$ is a polynomial function such that:
\begin{equation*}
F_{n+1} (0)=\dfrac{f^{(n+1)}(0)}{(n+1)!} \, ,\quad 
F_{n+1}' (0)=\dfrac{f^{(n+2)}(0)}{(n+2)!} \, .
\end{equation*}
Substituting the expression of $W_\flat$, and $W_\sharp$ in \eqref{ode-ref} 
shows (after a tedious computation!) that $(V_\flat,W_{\flat,n+1})$, and 
$(V_\sharp,W_{\sharp,n+1})$ are solutions to the following system of ODEs:
\begin{equation}
\label{eq:k}
\begin{cases}
V' = V \, (w_0+w_1 \, V +\cdots +w_n \, V^n +W_{n+1} \, V^{n+1}) \, ,& \\
W_{n+1}' =g_{n+1} (V,W_{n+1}) \, ,& 
\end{cases}
\end{equation}
where the function $g_{n+1}$ is given by:
\begin{align}
g_{n+1} (V,W_{n+1}) := &-(n+2) \, W_{n+1} \left( \sum_{k=0}^n w_k \, V^k 
+V^{n+1} \, W_{n+1} \right) -W_{n+1} \, \sum_{k=0}^n (k+1) \, w_k \, V^k \notag\\
&-W_{n+1} \, f(V) -F_{n+1} (V) \, \sum_{k=0}^n w_k \, V^k +b_n 
+\dfrac{f^{(n+1)}(0)}{(n+1)!} +V \, Q_{n+1}(V) \, ,\label{defgn+1}
\end{align}
and $Q_{n+1}$ is a polynomial function that satisfies:
\begin{equation*}
Q_{n+1}(0)=b_{n+1} +(n+4) \, w_1 \, w_{n+1} +f'(0) \, w_{n+1} 
+\dfrac{f^{(n+1)}(0)}{(n+1)!} \, w_1 +\dfrac{f^{(n+2)}(0)}{(n+2)!} \, w_0 \, .
\end{equation*}
When $n=1$, one has $Q_2 \equiv Q_2(0)=0$ (see the above definition for $b_2$). 
Using the expansion \eqref{expansion}, which is part of the induction assumption, 
we also know that:
\begin{equation*}
\lim_{\eta \rightarrow +\infty} (V_\flat,W_{\flat,n+1}) (\eta) 
=\lim_{\eta \rightarrow +\infty} (V_\sharp,W_{\sharp,n+1}) (\eta) 
=(0,w_{n+1}) =\Big( 0,\dfrac{b_n}{f(0)+(n+2)\, w_0} \Big) \, .
\end{equation*}

With the above definitions for $g_{n+1}$, and $Q_{n+1}$, we can check that 
$(0,w_{n+1})$ is a stationary solution to \eqref{eq:k}. (Recall that $w_{n+1}$ 
is defined as in Theorem \ref{smooth2} by the induction assumption.) We can 
also evaluate the Jacobian matrix of the vector field associated with the 
system of ODEs \eqref{eq:k}:
$$
\begin{pmatrix}
w_0 & 0 \\
b_{n+1} & -f(0)-(n+3) \, w_0 \end{pmatrix} =
\begin{pmatrix}
\lambda_1^{(n+1)} & 0 \\
b_{n+1} & \lambda_2^{(n+2)} \end{pmatrix} \, .
$$
There exists a positive number $a_{n+1} \le a_n$ such that for all $a \in 
(0,a_{n+1}]$, one has $\lambda_2^{(n+2)}<\lambda_1^{(n+2)}<0$, or equivalently 
$f(0)+(n+4) \, w_0>0$. In that case, the eigenvectors corresponding to the 
eigenvalues $\lambda_1^{(n+1)}$ and $\lambda_2^{(n+1)}$ are:
$$
e_1^{(n+1)} =\begin{pmatrix}
f(0)+(n+4) \, w_0 \\
b_{n+1} \end{pmatrix} \, ,\quad e_2^{(n+1)}=\begin{pmatrix}
0 \\
1 \end{pmatrix} \, .
$$
Using the same argument as in the proof of Proposition \ref{n=3}, we can 
conclude that the solutions $(V_\flat,W_{\flat,n+1})$, and 
$(V_\sharp,W_{\sharp,n+1})$ of \eqref{eq:k} are tangent to the straight line 
$(0,w_{n+1}) +\R \, e_1^{(n+1)}$ as $\eta$ tends to $+\infty$. In particular, 
this yields:
\begin{equation}
\label{lim0}
\lim_{\eta \rightarrow +\infty} \dfrac{W_{\flat,n+1}'(\eta)}{V_\flat'(\eta)} 
=\lim_{\eta \rightarrow +\infty} \dfrac{W_{\sharp,n+1}'(\eta)}{V_\sharp'(\eta)} 
=\dfrac{b_{n+1}}{f(0)+(n+4)\, w_0} =:w_{n+2} \, .
\end{equation}

\underline{Step 2}: let us define the function $\widetilde{w}_{n+1}$ by the formula:
\begin{equation*}
\widetilde{w}_{n+1} (\xi) :=\begin{cases}
W_{\flat,n+1} \circ \Xi_\flat^{-1} (\xi) &\text{if $\xi<0$,}\\
w_{n+1} &\text{if $\xi=0$,}\\
W_{\sharp,n+1} \circ \Xi_\sharp^{-1} (\xi) &\text{if $\xi>0$.}
\end{cases}
\end{equation*}
With this definition, $\widetilde{w}_{n+1}$ is continuous, and we have the relation:
\begin{equation}
\label{recurrence0}
w=w_0 +w_1 \, \hat{v} +\dots+w_n \, \hat{v}^n 
+\widetilde{w}_{n+1} \, \hat{v}^{n+1} \, .
\end{equation}
Moreover, using \eqref{lim0}, we obtain:
\begin{equation}
\label{recurrence4}
\lim_{\xi \rightarrow 0^-} \dfrac{\widetilde{w}_{n+1}' (\xi)}{\hat{v}'(\xi)}
=\lim_{\xi \rightarrow 0^+} \dfrac{\widetilde{w}_{n+1}' (\xi)}{\hat{v}'(\xi)}
=w_{n+2} \, ,
\end{equation}
which yields $\widetilde{w}_{n+1}'(0^+)=\widetilde{w}_{n+1}'(0^-)$. Therefore, we 
have $\widetilde{w}_{n+1} \in C^1 (\R)$. Moreover, using \eqref{eq:k}, we can compute:
\begin{equation}
\label{recurrence1}
\widetilde{w}_{n+1}' \, \hat{v} =g_{n+1} (\hat{v},\widetilde{w}_{n+1}) \, ,
\end{equation}
so we get $\widetilde{w}_{n+1}' \, \hat{v} \in C^1(\R)$.

\underline{Step 3}: we use an induction argument to show that $w \in C^{n+2} (\R)$ 
(which will imply immediately $\hat{v} \in C^{n+3} (\R)$). More precisely, we assume 
that for some $k \in \{ 0,\dots,n \}$, we have:
\begin{equation}
\label{recurrence2}
\widetilde{w}_{n+1} \, \hat{v}^k \in C^{k+1}(\R) \, ,\quad 
w \in C^{k+1} (\R) \, ,\quad 
\widetilde{w}_{n+1}' \, \hat{v}^{k+1} \in C^{k+1} (\R) \, .
\end{equation}
We are going to show that this property implies the same property with $k$ replaced 
by $k+1$. (Observe that step 2 above shows that the property \eqref{recurrence2} holds 
for $k=0$.)

We note that $\hat{v} \in C^{k+2} (\R)$, because $\hat{v}'=w \in C^{k+1} (\R)$. Moreover, 
we have $\widetilde{w}_{n+1} \, \hat{v}^{k+1} =(\widetilde{w}_{n+1} \, \hat{v}^k) \, 
\hat{v} \in C^{k+1} (\R)$, and we also have:
\begin{equation*}
(\widetilde{w}_{n+1} \, \hat{v}^{k+1})'=\widetilde{w}_{n+1}' \, \hat{v}^{k+1}
+(k+1) \, (\widetilde{w}_{n+1} \, \hat{v}^k) \, w \in C^{k+1} (\R) \, .
\end{equation*}
Therefore, we get $\widetilde{w}_{n+1} \, \hat{v}^{k+1} \in C^{k+2} (\R)$.

Using the relation \eqref{recurrence0}, we immediately obtain $w \in C^{k+2}(\R)$.

We have $\widetilde{w}_{n+1}' \, \hat{v}^{k+2}=(\widetilde{w}_{n+1}' \, \hat{v}^{k+1}) 
\, \hat{v} \in C^{k+1} (\R)$, and using \eqref{recurrence1}, we derive:
\begin{multline}
\label{recurrence3}
\big( \widetilde{w}_{n+1}' \, \hat{v}^{k+2} \big)' =\big( 
g_{n+1} (\hat{v},\widetilde{w}_{n+1}) \, \hat{v}^{k+1} \big)' \\
= (\partial_1 g_{n+1}) (\hat{v},\widetilde{w}_{n+1}) \, \hat{v}^{k+1} \, w 
+(\partial_2 g_{n+1}) (\hat{v},\widetilde{w}_{n+1}) \, \widetilde{w}_{n+1}' \, \hat{v}^{k+1} 
+(k+1) \, g_{n+1} (\hat{v},\widetilde{w}_{n+1}) \, \hat{v}^k \, w \, ,
\end{multline}
where $\partial_1 g_{n+1}$ (resp. $\partial_2 g_{n+1}$) denotes the partial derivative 
of $g_{n+1}$ with respect to its first (resp. second) variable. From the definition 
\eqref{defgn+1}, we see that $g_{n+1} (\hat{v},\widetilde{w}_{n+1})$ can be decomposed 
as follows:
\begin{equation*}
g_{n+1} (\hat{v},\widetilde{w}_{n+1})=-(n+2) \, \widetilde{w}_{n+1}^2 \, \hat{v}^{n+1} 
+\widetilde{w}_{n+1} \, P_1(\hat{v}) +P_0(\hat{v}) \, ,
\end{equation*}
where $P_0$, and $P_1$ are polynomial functions. Using this decomposition, and 
the induction assumption \eqref{recurrence2}, we can show that each term of the 
sum in the right-hand side of \eqref{recurrence3} belongs to $C^{k+1} (\R)$. 
Consequently $\widetilde{w}_{n+1}' \, \hat{v}^{k+2}$ belongs to $C^{k+2} (\R)$, 
and \eqref{recurrence2} holds with $k$ replaced by $k+1$. Because \eqref{recurrence2} 
holds for $k=0$, we get that \eqref{recurrence2} holds for $k=n+1$, so we have 
proved $w \in C^{n+2} (\R)$, and $\hat{v} \in C^{n+3}(\R)$.

\underline{Step 4}: it remains to show that $w$ satisfies the asymptotic expansion 
\eqref{expansion} at the order $n+1$. Using \eqref{recurrence4}, and 
$\widetilde{w}_{n+1} \in C^1 (\R)$, we obtain:
\begin{equation*}
\widetilde{w}_{n+1} (\xi)-w_{n+1}=w_{n+2} \, \hat{v} (\xi) +o(\hat{v} (\xi)) 
\, ,\quad \text{as $\xi \rightarrow 0$.}
\end{equation*}
Plugging this expansion in \eqref{recurrence0}, we obtain \eqref{expansion} at the 
order $n+1$, so the proof of the induction is complete.
\end{proof}

Once we know that the function $\hat{v}$ belongs to $C^{n+2} (\R)$, for 
$a \in (0,a_n]$, then $v=\hat{v}+(v_-+v_+)/2$ also belongs to $C^{n+2} (\R)$, 
and we have already seen in the previous section that $v$ does not vanish 
because $v(\xi)>v_+>0$ for all $\xi$. Moreover, the components $(\rho,u,e)$ 
of the shock profile are given by:
\begin{equation*}
\rho(\xi)=\dfrac{j}{v(\xi)} \, ,\quad u(\xi)=v(\xi)+\sigma \, ,\quad 
e(\xi)=\dfrac{(C_1-v(\xi))\, v(\xi)}{\gamma-1} \, ,
\end{equation*}
so one has $(\rho,u,e) \in C^{n+2} (\R)$, and the proof of Theorem \ref{smooth} 
is complete. (Recall that the strength of the shock tends to zero if, and only 
if $a=|u_+-u_-|/2$ tends to zero.)

\appendix
\section{Formal derivation of the model}
\label{model}

It is worth describing how the model \eqref{eulerevol}, \eqref{diffn} can be 
obtained from a more complete physical system. The derivation we propose below 
remains formal -- a rigorous proof being certainly delicate and beyond the scope 
of this work -- and we refer to \cite{BD, GL, LHM, MM} for further details. Let 
us introduce the specific intensity of radiation $f(t,x,v)$, that depends on 
a time variable $t\ge 0$, a space variable $x\in \R^N$, and a direction $v\in 
\SP^{N-1}$. We make the 'grey assumption', which means that the frequency 
dependence is ignored (all photons have the same frequency). Photons are 
subject to two main interaction phenomena:
\begin{itemize}
\item scattering produces changes in the direction of the photons,

\item absorption/emission where photons are lost/produced through a transfer 
mecanism with the surrounding gas.
\end{itemize}
The scattering phenomenon is described by the operator:
\[
Q_s(f)(t,x,v)=\sigma_s \, \Big( 
\displaystyle \int_{\SP^{N-1}} f(t,x,v') \, dv'-f(t,x,v) \Big) \, ,
\]
(with $dv$ the normalized Lebesgue measure on $\SP^{N-1}$), and the 
absorption/emission phenomenon is described by the operator:
\[
Q_a(f)(t,x,v)=\sigma_a \, \Big( 
\dfrac{\sigma}{\pi} \, \theta (t,x)^4 -f(t,x,v) \Big) \, ,
\]
where we used the Stefan-Boltzmann emission law, $\theta$ being the temperature 
of the gas, and $\sigma$ the Stefan-Boltzmann constant. In these definitions, 
the coefficients $\sigma_{s,a}$ are given positive quantities. These phenomena 
are both characterized by a typical mean free path, denoted $\ell_s,\ell_a$ 
respectively. Therefore, the evolution of the specific intensity is driven by:
\begin{equation}
\label{kindim}
\dfrac{1}{c} \, \partial_t f +v \cdot \nabla_x f =\dfrac{1}{\ell_s} \, Q_s(f) 
+\dfrac{1}{\ell_a} \, Q_a(f)=Q(f) \, ,
\end{equation}
where $c$ stands for the speed of light. The equation \eqref{kindim} is 
coupled to the Euler system describing the evolution of the fluid:
\begin{equation}
\label{eulerevolbis}
\begin{cases}
\partial_t \rho +\nabla_x \cdot (\rho \, u) =0 \, ,& \\
\partial_t (\rho \, u) +\nabla_x \cdot (\rho \, u\otimes u) +\nabla_x P
=-\dfrac{1}{c} \displaystyle \int_{\SP^{N-1}} v \, Q(f) \, dv \, ,& \\
\partial_t (\rho \, E) +\nabla_x \cdot (\rho \, E \, u+P\, u) 
=-\displaystyle \int_{\SP^{N-1}} Q(f) \, dv \, .&
\end{cases}
\end{equation}
The equations \eqref{kindim}, \eqref{eulerevolbis} are thus coupled by the 
exchanges of both momentum and energy, and by the Stefan-Boltzmann emission 
law. Observe that only the emission/absorption operator enters into the energy 
equation since the scattering operator is conservative (this would be different 
if Doppler corrections were taken into account). Note also that the total energy:
\[
\dfrac{1}{c} \displaystyle \int_{\R^N} \int_{\SP^{N-1}} f \, dv \, dx 
+\displaystyle \int_{\R^N} \rho \, E \, dx \, ,
\]
is (formally) conserved. Writing the system \eqref{eulerevolbis}, and the 
kinetic equation \eqref{kindim} in the dimensionless form, we can make four 
dimensionless parameters appear:
\begin{description}
\item - $\mathcal C$, the ratio of the speed of light over the typical sound 
speed of the gas,

\item - $\mathcal L_s$, the Knudsen number associated to the scattering,

\item - $\mathcal L_a$, the Knudsen number associated to the absorption/emission,

\item - $\mathcal P$, which compares the typical energy of radiation and the 
typical energy of the gas.
\end{description}
We thus obtain the rescaled equations:
\begin{equation}
\label{adim}
\begin{cases}
\dfrac{1}{\mathcal C} \, \partial_t f + v \cdot \nabla_x f 
=\dfrac{1}{\mathcal L_s} \, Q_s(f)+\dfrac{1}{\mathcal L_a} \, Q_a(f) \, ,& \\
\partial_t \rho +\nabla_x \cdot (\rho \, u)=0 \, ,& \\
\partial_t (\rho \, u) +\nabla_x \cdot (\rho \, u \otimes u) +\nabla_x P 
=\dfrac{\mathcal P}{\mathcal L_s} \, \sigma_s \, 
\displaystyle \int_{\SP^{N-1}} v \, f(v) \, dv \, ,& \\
\partial_t (\rho \, E) +\nabla_x \cdot (\rho \, E \, u+P\, u) 
=-\dfrac{\mathcal P}{\mathcal L_a} \, \sigma_a \, \Big( 
\theta^4 -\displaystyle \int_{\SP^{N-1}} f (v) \, dv \Big) \, .&
\end{cases}
\end{equation}
System \eqref{eulerevol}, \eqref{diffn} is then obtained in two steps. First of 
all, we assume $\mathcal C \gg 1$. Next, we keep $\mathcal P$ of order 1, and we 
are concerned here with a regime where scattering is the leading phenomenon: the 
mean free paths are rescaled according to:
\[
\mathcal L_s \simeq \dfrac{1}{\mathcal C} \, ,\qquad 
\mathcal L_a \simeq \mathcal C \, .
\]
The asymptotics can be readily understood by means of the Hilbert expansion:
\[
f=f^{(0)} +\dfrac{1}{\mathcal C} \, f^{(1)} +\dfrac{1}{\mathcal C^2} \, f^{(2)}+\dots
\]
Identifying the terms arising with the same power of $1/\mathcal C$, we get:
\begin{description}
\item - at the leading order, $f^{(0)}$ belongs to the kernel of the scattering 
operator, so that is does not depend on the microscopic variable $v$: 
$f^{(0)}(t,x,v)=n(t,x)$,

\item - the relation $Q_s(f^{(1)})=v \cdot \nabla_x f^{(0)}$ then leads to: 
$f^{(1)}(t,x,v)=-\frac{1}{\sigma_s}\, v \cdot \nabla_x n(t,x)$,

\item - integrating the equation for $f^{(2)}$ over the sphere yields:
\[
\partial_t n -\dfrac{1}{N \, \sigma_s} \, \Delta_x n 
=\sigma_a \, ( \theta^4 -n) \, .
\]
\end{description}
Note also that in the momentum equation, we have:
\[
\dfrac{\sigma_s}{\mathcal L_s} \displaystyle \int_{\SP^{N-1}} v \, f(v) \, dv 
\simeq \sigma_s \, \displaystyle \int_{\SP^{N-1}} v \, f^{(1)}(v) \, dv 
=-\dfrac{1}{N} \, \nabla_x n \, .
\]
Finally, we obtain the limit system:
\begin{equation}
\label{lim1}
\begin{cases}
\partial_t \rho +\nabla_x \cdot (\rho \, u)=0\, ,& \\
\partial_t (\rho \, u) +\nabla_x \cdot (\rho \, u \otimes u) +\nabla_x P 
=-\dfrac{\mathcal P}{N} \, \nabla_x n \, ,& \\
\partial_t (\rho \, E) +\nabla_x \cdot (\rho \, E \, u +P\, u) 
=-\mathcal P \, \sigma_a \, ( \theta^4-n ) \, ,& \\
\partial_t n -\dfrac{1}{N \, \sigma_s} \, \Delta_x n =\sigma_a \, 
( \theta^4-n ) \, .
\end{cases}
\end{equation}
The system \eqref{lim1} describes a nonequilibrium regime, where the material 
and the radiations have different temperatures ($\theta\neq n^{1/4}$); the 
equilibrium regime would correspond to assuming that the emission/absorption 
is the leading contribution.

After this first asymptotics, we perform a second asymptotics where we set:
\[
\mathcal P \ll 1 \, ,\qquad \mathcal P \, \sigma_a=1 \, ,\qquad 
N\, \sigma_s=1/\sigma_a \, .
\]
This leads to \eqref{eulerevol}, \eqref{diffn}. Of course, one might wonder 
how this second approximation modifies the shock profiles compared to 
\eqref{lim1}, in particular when we get rid of the radiative pressure 
in the momentum equation. We refer to \cite[page 579]{MM} for some aspects 
of this problem.

\bibliographystyle{alpha}
\bibliography{lcg}
\end{document}